\documentclass[10
pt,reqno,fleqn]{amsart}
\usepackage{amsfonts,amsmath,amsthm,amssymb}
\usepackage{pdfsync}

\usepackage{mathrsfs}

\newtheorem{thm}{Theorem}

\newtheorem{lemma}{Lemma}

\newcommand{\R}{\ensuremath{\mathbb{R}}}
\newcommand{\Z}{\ensuremath{\mathbb{Z}}}
\newcommand{\C}{\ensuremath{\mathbb{C}}}

\newcommand{\1}{\ensuremath{\mathbf{1}}}

\newcommand{\<}{\ensuremath{\lesssim}}

\newcommand{\eps}{\ensuremath{\varepsilon}}
\newcommand{\la}{\ensuremath{\lambda}}
\newcommand{\La}{\ensuremath{\Lambda}}

\newcommand{\eq}{\begin{equation}}
\newcommand{\ee}{\end{equation}}

\numberwithin{equation}{section}

\begin{document}
\title{Discrete multilinear spherical averages}
\author{Brian Cook}

\begin{abstract}
In this note we give a characterization of $\ell^{p}\times ...\times \ell^{p}\to\ell^q$ boundedness of maximal operators associated to multilinear convolution averages over spheres in $\Z^n$.\end{abstract}

\thanks{
The author was supported in part by NSF grant DMS1147523}

\maketitle

\section{Introduction}

Consider the multilinear convolution operators which are defined by \[
A_\la[f_1,...,f_n](y)=\frac{1}{r(\la)}\sum_{|x|^2=\la} f_1(y-x_1)...f_n(y-x_n)\]
where $r(\la)=\{x\in\Z^n: |x|^2=x_1^2+x_2^2+...+x_n^2\}$ and $f_1,...,f_n$ are functions defined on the integers. The normalization factor $r(\la)$ is well known to be well behaved when $n\geq5$. In particular, for these $n$ there are constants $c_n$ and $C_n$ such that  $c_n\la^{n/2-1}\leq r(\la) \leq C_n\la^{n/2-1}$. When $n=4$ this type of regularity disappears, but still each $A_\la$ is well defined. If $n<4$ the operators are not defined for all $\la\geq1$ and the restriction $n>3$ is assumed throughout. 

The operators $A_\la$ are the discrete analogue of the  multilinear convolution operators defined in the Euclidean setting considered in \cite{Dan}.  In that work the boundedness on $L^p(\R)\times...\times L^p(\R)\to L^q(\R)$ was characterized for a single radius $\la$. This particular question in our context is not overly interesting. Indeed, a later paper of Oberlin, \cite{Dan2}, covers a multilinear version of Young's inequality that can be directly applied to address this question for the $A_\la$ given above. A more interesting problem is to consider $\ell^p(\Z)\times...\times \ell^p(\Z)\to \ell^q(\Z)$ boundedness of the maximal operators defined pointwise by \[
A_*[f_1,...,f_n](y)=\sup_{\la\geq 1}|A_\la[f_1,...,f_n](y)|,\]
which is the motivation of this note.

From a slightly different point of view we observe that the operators $A_\la$ are the multilinear analogues of the discrete spherical averages considered in \cite{MSW} that are given by \[
S_\la \phi(y)=\frac{1}{r(\la)}\sum_{|x|^2=\la} \phi(y-x),\]
where $\phi:\Z^n\to\C$. We have \[
A_\la[f_1,...,f_n](y)  =S_\la \Phi(\tilde{y})\]
where $ \Phi =f_1\otimes...\otimes f_n$ and $\tilde{y}=(y,...,y)$ is the image of $y$ under the natural embedding of $\Z$ into the diagonal of $\Z^n$. The associated maximal operators $S_*$ are known to be bounded on $\ell^p(\Z)$ if and only if $n\geq 5$ and $p> n/(n-2)$. This point of view does not yield any immediate results for us, but it does  motivate the notion that one should be able to multilinearize  the methods used to study $S_*$ to obtain results for the maximal operators considered here. This is probably the right way to go in the end, due to certain applications, but we are not sure at this point if the required Euclidean analogues needed to carry this out are currently known. Here we consider a much simpler approach which does not rely on any results in the $\R^n$ setting.  In fact, we give a characterization of $\ell^{p}\times ...\times \ell^{p}\to\ell^q$ boundedness of the operators $A_*$ via a reduction to the discrete Hardy-Littlewood maximal theorem. The form of this maximal theorem we employ is as follows. 

\bigskip
\begin{lemma}\label{lem1}
Let $f$ be a function defined on $\Z$. We have the inequality\footnote{The notation $f\<g$ means that there is a constant $C$ such that $|f|\leq C g$ where $g\geq0$. The constants $C$ are allowed to depend on any parameters other than than those specially related $\la$ or $N$.}\[
||\sup_{N\geq 1} |N^{-1}\sum_{x\in[-N,N]} f(y-x)|\,||_{\ell^p}\< ||f||_{\ell^p}\]
for all $p>1$. 
\end{lemma}
\bigskip

This will provide us with our desired result.

\bigskip
\begin{thm}\label{thm1}
Let $1\leq q <\infty$. The operator $A_*$ is bounded from $\ell^{p}\times ... \times \ell^{p} $ to $\ell^q$ if and only if $n\geq 5$ and $1\leq p\leq nq$. 
\end{thm} 
\bigskip

The reduction of Theorem \ref{thm1} to the discrete Hardy-Littlewood maximal inequality is an application of a type of discrete restriction inequality. The specific result we are interested in is provided in \cite{Bourgain}.

\bigskip
\begin{lemma}\label{lem2}
Let $f$ be a function defined on $\Z$ and $n> 4$. Then\[
\int_{\Pi} \left|\sum_{x\in[-N,N]}f(x)e(\alpha x^2)\right|^n\,d\alpha\< N^{n-2}\left(N^{-1} \sum_{x\in[-N,N]}|f(x)|^2\right)^{n/2}.\]
\end{lemma}

\bigskip

Lemma \ref{lem2} has higher degree counterparts, meaning the $x^2$ in the phase is replaced by $x^d$ for some $d>2$. Results of this type are addressed in a recent paper of K. Hughes and K. Henriot \cite{k2}. Further results are also known, for example the paper \cite{Wooley} by T. Wooley provides restriction type results related to Vinogradov's mean value theorem.  Such results can be applied to obtain analogues of Theorem \ref{thm1}, although we do not consider such things here.

\bigskip
\section{Proof}

We first dispense with the necessary conditions. The first simple observation is that if $A_*$ is bounded from $\ell^{p}\times ... \times \ell^{p} $ to $\ell^q$, then it is bounded from $\ell^{p'}\times ... \times \ell^{p'} $ to $\ell^q$ when $p'<p$. 
 
Next we consider the case when the functions $f_1=...=f_n=\phi$, where $\phi$ is the unit mass at the origin, i.e. \[
\phi(x) = \left\{
        \begin{array}{ll}
            1 &  x = 0 \\
            0 & x \not= 0\\
        \end{array}.
    \right.\] 
In this case we have \[
A_*[\phi,...,\phi](y)=\sup_{\la\geq1}\,\frac{1}{r(\la)}\sum_{|x|^2=\la}\phi(y-x_1)...\phi(y-x_n)=\frac{1}{r(ny^2)}.\]
If $n=4$ then we have that $r(4^k)=24$, implying that  $A_*[\phi,...,\phi](y)=1/24$ when $y$ is a power of two. Hence $A_*[\phi,...,\phi]$ does not belong to any $\ell^q(\Z)$ space when $q<\infty$.

The final requirement, namely that  $p\leq qn $,  is more or less a scaling issue. To see precisely why this condition is required we consider the scenario when $f_1=...=f_n=\chi:=\1_{[-2M,2M]}$ for a fixed large integer $M$.  For $|y|\leq M$\[
A_*[\chi,...,\chi](y)=\sup_{\la\geq1}\,\frac{1}{r(\la)}\sum_{|x|^2=\la}\chi(y-x_1)...\chi(y-x_n),\]
 is at least \[
\frac{1}{r(M^2)}\sum_{|x|^2=M^2}\chi(y-x_1)...\chi(y-x_n)=1.\]
Then we have that \[
||A_*[\chi,...,\chi]||_{\ell^q(\Z)}\geq M^{1/q}.\]
On the other hand we have that\[ 
||\chi||_{\ell^p(\Z)}^n\<M^{n/p}.\]
The conclusion of Theorem \ref{thm1} is then seen to be false when $p>nq$ by selecting $M$ sufficiently large. 

We now proceed with the proof in the remaining case when $p=nq$ and $n\geq5$.  This begins with the standard observation that we can write \[
A_\la[f_1,...,f_n](y)=\frac{1}{r(\la)}\int_\Pi\sum_{|x_1|\leq N_\la}...\sum_{|x_n|\leq N_\la}f_1(y-x_1)...f_n(y-x_n)e((|x|^2-\la)\alpha)\,d\alpha\]
where  $N_\la$ is chosen to be $\approx \la^{1/2}$ (for example, one can choose $N_\la$ to be the closest integer to $2\la^{1/2}$), $\Pi=\R/\Z$, and $e(x)=e^{2\pi i x}$. Defining the exponential sums \[
W_{i,\la}(\alpha,y)=\sum_{|x_i|\leq N_\la}f_i(y-x_i)e(x_i^2\alpha)\]
puts this in the form \[
\frac{1}{r(\la)}\int_\Pi   W_{1,\la}(\alpha,y)...W_{n,\la}(\alpha,y)e(-\la)\alpha)\,d\alpha.\]

We can now apply H\"older's inequality to get the bound \[
|A_\la[f_1,...,f_n](y)|\leq\frac{1}{r(\la)}\prod_{i=1}^n\left( \int_\Pi |W_{i,\la}(\alpha,y)|^n\,d\alpha\right)^{1/n}.\]
This puts us in a position to apply Lemma \ref{lem2} which results in\[
|A_\la[f_1,...,f_n](y)|\leq\frac{N_\la^{n-2}}{r(\la)}\prod_{i=1}^n\left(N_\la^{-1}\sum_{|x_i|\leq N_\la} |f_i(y-x_i)|^2\right)^{1/2}.\]
Notice that \[
\sup_{\la\geq 1}\frac{N_\la^{n-2}}{r(\la)}\] is uniformly bounded in $N$,
which gives that \[
|A_*[f_1,...,f_n](y)|^q\<\prod_{i=1}^n\left(\sup_{\la\geq1}\,N_\la^{-1}\sum_{|x_i|\leq N_\la} |f_i(y-x_i)|^2\right)^{q/2},\]
or just \[
|A_*[f_1,...,f_n](y)|^\<\prod_{i=1}^n\left(\sup_{N\geq1}\,N^{-1}\sum_{|x_i|\leq N} |f_i(y-x_i)|^2\right)^{q/2}.\]

We proceed by summing  in $y$. This gives\[
||A_*[f_1,...,f_n]||_{\ell^q}^q\< \sum_{y\in \Z} \,\prod_{i=1}^n\left(\sup_{N\geq1}\,N^{-1}\sum_{|x_i|\leq N} |f_i(y-x_i)|^2\right)^{q/2}.\]
Another application of H\"older's inequality then gives \[
||A_*[f_1,...,f_n]||_{\ell^q}^q\< \, \prod_{i=1}^n\left(\sum_{y\in \Z}\left(\sup_{N\geq1}\,N^{-1}\sum_{|x_i|\leq N} |f_i(y-x_i)|^2\right)^{(qn/2)}\right)^{1/n}.\]

Recall now that $q=pn$, and then the terms \[
\sum_{y\in \Z}\left(\sup_{N\geq1}\,N^{-1}\sum_{|x_i|\leq N} |f_i(y-x_i)|^2\right)^{qn/2}\]
are at most \[
\sum_{x\in\Z} |f_i(x)|^{p}=||f_i||^{p}_{\ell^{p}}\]
by the Hardy-Littlewood maximal inequality applied with the functions $|f_i|^2$. In turn we have \[
||A_*[f_1,...,f_n]||^q_{\ell^q}\<\prod_{i=1}^n||f_i||^{(p/n)}_{\ell^{p}}=\prod_{i=1}^n||f_i||^{q}_{\ell^{p}},\]
which is the desired result.

\vskip0.2in
\noindent \author{\textsc{Brian Cook}}\\
Department of Mathematical Sciences \\
Kent State University\\
Kent, OH, USA\\
Electronic address: \texttt{briancookmath@gmail.com}

\end{document}